\documentclass{article} 
\usepackage{graphicx}
\usepackage{layout}
\usepackage{fullpage}
\usepackage{amsmath}
\usepackage{amssymb}
\usepackage{float}
\usepackage{hyperref}
\title{A Generalization of Siegel’s Method to Jacobi's $\vartheta_1$ Function}
\date{}
\author{Maher Mamah $\&$ Ali Saraeb}
\begin{document}
\maketitle
\begin{figure}[H]
    \centering
    \includegraphics[scale=0.35]{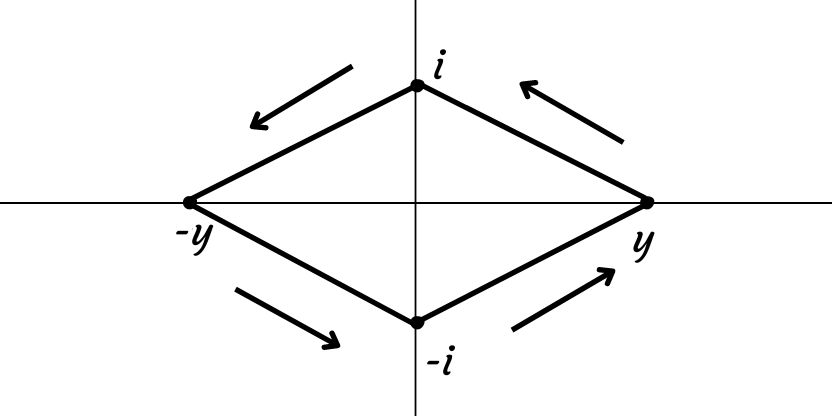}
    \caption{Contour \( C \)}
    \label{fig 1}
\end{figure}

\textbf{Abstract.} We present a new proof of the transformation law of $\vartheta_1$ under the action of the generator of the full modular group $\Gamma$ using Siegel's method. \\ 
\begin{center}
    1. INTRODUCTION
\end{center}
Historically, the study of elliptic functions has taken many different directions, stemming from the great fathers of this theory like Abel, Weierstrass, Jacobi etc. However, Jacobi theta functions form an important study subject in this theory, for they are the building blocks of elliptic functions.\\ The most substantial study of these theta function is in the realm of the theory of modular forms, where these function are studied under the action of the elements of the modular group $\Gamma$ and its subgroups. And as the theory evolved with other branches of analysis, newer and lighter proofs were occasionally introduced as new ways to approach the laws of transformation of these theta functions.\\
Define $\mathbb{H}=\{x+iy: y>0\}$ to be the upper half-plane.  For $\tau\in \mathbb{H}, z\in \mathbb{C}$, let $q=e^{\pi i\tau}$ and $w=e^{\pi iz}$, and we define 
\begin{align}\label{eq:1}
\vartheta_1(z,\tau) = -i w q^{1/4} \prod_{n=1}^{\infty} (1-q^{2n})(1- w^2 q^{2n})(1- w^{-2} q^{2n-2}).
 \end{align} \\
The transformation law of $\vartheta_1$  under the action of the inversion matrix $S=\left(\begin{array}{cc}
   0  &  -1\\
    1 & 0
\end{array}\right)$ is given by
\begin{align} \label{eq:2}
\displaystyle
\vartheta_1\left(\frac{z}{\tau},\frac{-1}{\tau}\right)= -i (-i\tau)^{1/2} e^{\pi i \frac{z^2}{\tau}} \vartheta_1(z,\tau),
\end{align} \\
and our main goal is to give a new proof of this property.\\
\\
Using residue calculus, Siegel presented the transformation law of the Dedekind eta function under the action of the inversion matrix $ \cite{Siegel's method}\relax$. Several extensions of Siegel's proof have been introduced, such as Raji's extension \cite{Raji} and the extension by Me'Meh and Saraeb \cite{DefTheta}. Inspired by Siegel, Raji proved in $ \cite{Raji}\relax \text{ }$the transformation law for another one of the Jacobi theta functions, $\vartheta_3$, under the inversion matrix. In this paper, we  are complementing Raji's paper by doing the same for $\vartheta_1$. For the remaining theta functions $\vartheta_2$ and $\vartheta_4$ they are related to $\vartheta_1$ and $\vartheta_3$ via the following relations described in [\cite{Rademacher}\relax- Chapter X]
\[\vartheta_3\left(z+\frac{1}{2},\tau\right)=\vartheta_4(z,\tau)~~~~~~~\vartheta_2\left(z+\frac{1}{2},\tau\right)=-\vartheta_1(z,\tau),\]
so by applying a sequence of changes of variables we get similar results for $\vartheta_2$ and $\vartheta_4$. It is important to note that Raji in [3] proved the behavior of $\vartheta_3$ under the inversion since the behavior under the translation matrix is clear from the Fourier expansion. Similarly, the modularity of $\vartheta_1$  with respect to the full group $SL_2(\mathbb{Z})$ will follow from (\ref{eq:2}) and its Fourier Expansion. 
\begin{center}
  2. THE TRANSFORMATION LAW FOR $\vartheta_1$
\end{center}
\textbf{Theorem 1.}\label{theo1} For $\tau\in \mathbb{H}, z\in \mathbb{C}$, the transformation formula for $\vartheta_1$ is the one presented in equation (2), namely
\begin{align*}
\displaystyle
\vartheta_1(\frac{z}{\tau},\frac{-1}{\tau})= -i (-i\tau)^{1/2} e^{\pi i \frac{z^2}{\tau}} \vartheta_1(z,\tau) .
\end{align*} \\[1\baselineskip]
If $z \in \mathbb{C} \backslash \{0\}$ and $x \in \mathbb{C}$, we set $z^x = e^{x \text{log} z}$ where $- \pi <  \text{arg} (z) \leq \pi$. \\  \\
\textbf{Definition 1.} For $\tau\in \mathbb{H} \text{ and } z\in \mathbb{C}$, let \( \phi(z,\tau) = \log(\vartheta_1(z,\tau)) - \log\left(\vartheta_1\left(\frac{z}{\tau}, \frac{-1}{\tau}\right)\right) \). \\ \\
\textbf{Remark 1.} $\phi$ takes the value infinity when $\vartheta_1$ vanishes, namely when $\displaystyle z= m +n \tau$ for $m, n \in \mathbb{Z}$\\ \\
\textbf{Lemma 1.\label{lemma 1}} Let  $b<0, ~ 0 <a<1$,~ $y>|b|,$ and $z= a+ib$. Then we have 
\begin{align*}
\displaystyle \phi(z, iy) &= 
\sum_{m=1}^{\infty} \frac{1}{m} \frac{1}{1-e^{2m\pi y}} + \sum_{m=1}^{\infty} \frac{e^{2m\pi i z}}{m} \frac{1}{1-e^{2m\pi y}}  + \sum_{m=1}^{\infty} \frac{e^{-2m\pi i z} }{m} \frac{e^{2m\pi y}}{1-e^{2m\pi y}}  \nonumber \\ &- \sum_{m=1}^{\infty} \frac{1}{m} \frac{1}{1-e^{2m\pi \frac{1}{y} }} - \sum_{m=1}^{\infty} \frac{e^{2m\pi \frac{z}{y}}}{m} \frac{1}{1-e^{2m\pi \frac{1}{y}}} \nonumber - \sum_{m=1}^{\infty} \frac{e^{-2m\pi \frac{z}{y}} }{m} \frac{e^{2m\pi \frac{1}{y}}}{1-e^{2m\pi \frac{1}{y}}}  \\ &- \frac{\pi z}{v}  + \pi i z -\frac{\pi }{4}(y - \frac{1}{y}) .\end{align*} 
\textit{Proof.} \\ \\
We first observe that the logarithm of Equation (\ref{eq:1}) gives\\
\begin{align} 
\text{log}(\vartheta_1(z,\tau))= \text{log}(-i w q^{1/4})+ \sum_{n=1}^{\infty} \text{log}(1-q^{2n}) + \sum_{n=1}^{\infty} \text{log}(1- w^2 q^{2n}) + \sum_{n=1}^{\infty} \text{log}(1- w^{-2} q^{2n-2}).
\end{align}
If $\tau=iy$, we get $ |q^{2n}|= e^{-2 \pi n y} ,\text{ } |w^2 q^{2n}| = e^{-2 \pi (n y +  b) }, \text{ }$  $|w^{-2} q^{2n-2}|= e^{-2 \pi ((n-1) y -  b) } \text{ }$ for all $n \geq 1$, then under the assumptions $b<0$, and $y>|b|$, we obtain the condition that all the latter quantities are strictly less than 1.
Now, we expand the form of $\text{log}(\vartheta_1(z,\tau))$  in (3) by using the Taylor expansion of log$(1-t)$ for $|t|<1$ as follows.\\
\begin{align*}
\text{log}(\vartheta_1(z,iy))& = \text{log}(-i w q^{1/4})- \sum_{n=1}^{\infty} \sum_{m=1}^{\infty} \frac{1}{m}(q^{2n})^m - \sum_{n=1}^{\infty} \sum_{m=1}^{\infty} \frac{1}{m} (w^2 q^{2n})^m  - \sum_{n=1}^{\infty} \sum_{m=1}^{\infty} \frac{1}{m} (w^{-2} q^{2n-2})^m \nonumber \\
&= \text{log}(-i w q^{1/4}) - \sum_{m=1}^{\infty} \frac{1}{m} \frac{q^{2m}}{1-q^{2m}} - \sum_{m=1}^{\infty} \frac{w^{2m}}{m} \frac{q^{2m}}{1-q^{2m}} - \sum_{m=1}^{\infty} \frac{w^{-2m} q^{-2m}}{m} \frac{q^{2m}}{1-q^{2m}} \nonumber \\
&= \text{log}(-i e^{\pi i z} e^{\frac{-\pi y}{4}}) - \sum_{m=1}^{\infty} \frac{1}{m} \frac{e^{-2m\pi y}}{1-e^{-2m\pi y}} - \sum_{m=1}^{\infty} \frac{e^{2m\pi i z}}{m} \frac{e^{-2m\pi y}}{1-e^{-2m\pi y}} \nonumber \\ 
&- \sum_{m=1}^{\infty} \frac{e^{-2m\pi i z} }{m} \frac{1}{1-e^{-2m\pi y}} \nonumber\\
&= \text{log}(-i e^{\pi i z} e^{\frac{-\pi y}{4}}) + \sum_{m=1}^{\infty} \frac{1}{m} \frac{1}{1-e^{2m\pi y}} + \sum_{m=1}^{\infty} \frac{e^{2m\pi i z}}{m} \frac{1}{1-e^{2m\pi y}} \nonumber \\ 
&+ \sum_{m=1}^{\infty} \frac{e^{-2m\pi i z} }{m} \frac{e^{2m\pi y}}{1-e^{2m\pi y}}.
\end{align*}  
Similarly, if we work under the conditions $y>|b|$ (particularly $y>0$) and $0<a<1$, we get that all of the real numbers $ \displaystyle e^{-2 \pi \frac{n}{y}} ,\text{ } e^{-2 \pi \frac{-a+n}{y}  }, \text{ and}$  $\displaystyle e^{-2 \pi \frac{a+n-1}{y}  } \text{ }$ are strictly less than 1 for all $n \geq 1$. Thus, we can ana\text{log}ously determine the form of $ 
\displaystyle \text{log}(\vartheta_1(\frac{z}{iy},\frac{i}{y}))$ from (4) as follows 
\begin{align}
\text{log}\left(\vartheta_1\left(\frac{z}{iy},\frac{i}{y}\right)\right)&= \text{log}\left(-i e^{\pi \frac{z}{y}} e^{\frac{-\pi }{4y}}\right) + \sum_{m=1}^{\infty} \frac{1}{m} \frac{1}{1-e^{2m\pi \frac{1}{y} }} + \sum_{m=1}^{\infty} \frac{e^{2m\pi \frac{z}{y}}}{m} \frac{1}{1-e^{2m\pi \frac{1}{y}}} \nonumber \\ 
&+ \sum_{m=1}^{\infty} \frac{e^{-2m\pi \frac{z}{y}} }{m} \frac{e^{2m\pi \frac{1}{y}}}{1-e^{2m\pi \frac{1}{y}}} .
\end{align} 
Subtracting (4) from (5), we obtain 
\begin{align*}
\displaystyle \phi(z, iy) &= \text{log}(\vartheta_1(z,iy)) -\text{log}\left(\vartheta_1\left(\frac{z}{iy},\frac{i}{y}\right)\right) \\ &= 
 \sum_{m=1}^{\infty} \frac{1}{m} \frac{1}{1-e^{2m\pi y}} + \sum_{m=1}^{\infty} \frac{e^{2m\pi i z}}{m} \frac{1}{1-e^{2m\pi y}}  + \sum_{m=1}^{\infty} \frac{e^{-2m\pi i z} }{m} \frac{e^{2m\pi y}}{1-e^{2m\pi y}}  \nonumber \\ &- \sum_{m=1}^{\infty} \frac{1}{m} \frac{1}{1-e^{2m\pi \frac{1}{y} }} - \sum_{m=1}^{\infty} \frac{e^{2m\pi \frac{z}{y}}}{m} \frac{1}{1-e^{2m\pi \frac{1}{y}}} \nonumber - \sum_{m=1}^{\infty} \frac{e^{-2m\pi \frac{z}{y}} }{m} \frac{e^{2m\pi \frac{1}{y}}}{1-e^{2m\pi \frac{1}{y}}}  \\ &+ \text{log}(-i e^{\pi i z} e^{\frac{-\pi y}{4}})   -\text{log}(-i e^{\pi \frac{z}{y}} e^{\frac{-\pi }{4y}})
 \\ &= 
\sum_{m=1}^{\infty} \frac{1}{m} \frac{1}{1-e^{2m\pi y}} + \sum_{m=1}^{\infty} \frac{e^{2m\pi i z}}{m} \frac{1}{1-e^{2m\pi y}}  + \sum_{m=1}^{\infty} \frac{e^{-2m\pi i z} }{m} \frac{e^{2m\pi y}}{1-e^{2m\pi y}}  \nonumber \\ &- \sum_{m=1}^{\infty} \frac{1}{m} \frac{1}{1-e^{2m\pi \frac{1}{y} }} - \sum_{m=1}^{\infty} \frac{e^{2m\pi \frac{z}{y}}}{m} \frac{1}{1-e^{2m\pi \frac{1}{y}}} \nonumber - \sum_{m=1}^{\infty} \frac{e^{-2m\pi \frac{z}{y}} }{m} \frac{e^{2m\pi \frac{1}{y}}}{1-e^{2m\pi \frac{1}{y}}}  \\ &- \frac{\pi z}{y}  + \pi i z -\frac{\pi }{4}\left(y - \frac{1}{y}\right).\end{align*} 
This completes the proof of \textit{Lemma 1}.$~~~~~~~~~~~~~~~~~~~~~~~~~~~~~~~~~~~~~~~~~~~~~~~~~~~~~~~~~~~~~~~~~~~~~~~~~~~~~~~~~~~~~~~~~~~~~~~~~~\blacksquare$ \\ \\
\textbf{Remark 2.} The conditions $y>|b|$ with $b<0<a<1$, where $z= a+ib$ guarantee the convergence of the latter series. \\ \\
\textbf{Definition 2.} Let $n$ be a positive integer, set $N= n+ \frac{1}{2}$, $y>|b|$ with $b<0<a<1$, and $z= a+ib$, we define a meromorphic function $F_n$ as
\begin{align} F_n(\zeta)= -\frac{1}{8\zeta} \text{cot}\left(\pi i N \zeta\right) \text{cot}\left(\pi N \frac{1}{y}\zeta\right) +\frac{1}{\zeta}\left(\frac{1}{1-e^{2\pi N \zeta}}\right)\left(\frac{e^{-2\pi i (-\frac{zN}{y}+\frac{N}{y})\zeta}}{1-e^{-2\pi i \frac{N}{y}\zeta}}\right), \end{align} 
which exhibits a pole of order 3 at 0 and simple poles at $\displaystyle \displaystyle \frac{ik}{N}$ and $\displaystyle \frac{ky}{N}$ for k $\in~\mathbb{Z}\backslash\lbrace 0\rbrace$. \\ \\
\indent
Now, we define the contour \( C \) shown in Figure \ref{fig 1} for the residue theorem application. First, we analyze the residues enclosed by \( C \). 

\noindent
\textbf{Lemma 2.} The sum of the residues of \( F_n \) at its poles within the contour \( C \), multiplied by \( 2 \pi i \), yields the following expression:
\begin{align*}
 2\pi i \cdot \text{Res}_{C}[F_n] & = \sum_{k=1}^{n} \frac{1}{k} \frac{1}{1-e^{2k\pi y}} + \sum_{k=1}^{n} \frac{e^{2k\pi i z}}{k} \frac{1}{1-e^{2k\pi y}}  \quad + \sum_{k=1}^{n} \frac{e^{-2k\pi i z} }{k} \frac{e^{2k\pi y}}{1-e^{2k\pi y}} \\
& \quad - \sum_{k=1}^{n} \frac{1}{k} \frac{1}{1-e^{2k\pi \frac{1}{y} }} - \sum_{k=1}^{n} \frac{e^{2k\pi \frac{z}{y}}}{k} \frac{1}{1-e^{2k\pi \frac{1}{y}}}  \quad - \sum_{k=1}^{n} \frac{e^{-2k\pi \frac{z}{y}} }{k} \frac{e^{2k\pi \frac{1}{y}}}{1-e^{2k\pi \frac{1}{y}}} \\
& \quad - \frac{\pi}{4} \left(y -\frac{1}{y}\right) + \pi i z + \frac{\pi z^2}{y} - \frac{\pi z}{y} - \frac{\pi i}{2}.
\end{align*}
where \( \phi \) was defined in \textit{Definition 1}, and
\[ \text{Res}_{C}[F_n] = \text{Res}[ F_n , 0 ] + \sum_{\substack{k=-n \\ k \ne 0}}^{n} \text{Res}[ F_n, \frac{ik}{N} ] + \text{Res}[ F_n, \frac{ky}{N} ]. \]
\textbf{Remark 3.} To gain insight into the direction of this analysis, we encourage the reader to compare the results of \textit{Lemma 1} and \textit{Lemma 2}. \\ \\
\textit{Proof of Lemma 2.}\\
We start by computing the residues of $F_n$ at the poles $\zeta=0$, $\zeta= \displaystyle \frac{ik}{N}$, and $ \zeta=\displaystyle \frac{ky}{N} $ for k $\in\mathbb{Z}\backslash\lbrace 0\rbrace$. Also note that the properties of the exponential and the cotangent imply that $F_n$ has no poles except at the asserted points. \\ \\
Using the Laurent expansions of $\displaystyle \text{cot}(\zeta)$ and $\displaystyle \frac{1}{1-e^\zeta}$ around zero, namely \\ \\
\[\displaystyle \text{cot}(\zeta) = \frac{1}{\zeta} - \frac{\zeta}{3} +O(\zeta^3)\]
and
\[\displaystyle \frac{1}{1-e^\zeta} = -\frac{1}{\zeta} +\frac{1}{2} -\frac{\zeta}{12} +O(\zeta^3),\]\\ \\
we obtain the Laurent expansions of the following terms \\ \\ 
 \[\displaystyle -\frac{1}{8\zeta} \text{cot}(\pi i N \zeta) \text{cot}\left(\pi N \frac{1}{y}\zeta\right) = \frac{iy}{8 \pi^2 N^2} \frac{1}{\zeta^3} + \frac{i}{24}\left(y -\frac{1}{y}\right) \frac{1}{\zeta} +O( \zeta ) ~~~~~~~~~~~(R_1) \] 
\\
\[ \frac{1}{\zeta}\left(\frac{1}{1-e^{2\pi N \zeta}}\right)\left(\frac{e^{-2\pi i (-\frac{zN}{y}+\frac{N}{y})\zeta}}{1-e^{-2\pi i \frac{N}{y}\zeta}}\right) = \frac{iy}{4 \pi^2 n^2} \frac{1}{\zeta^3} - \frac{1}{4 \pi N} (iy +2z-1 ) \frac{1}{\zeta^2} + \left(\frac{i}{12}\left(y- \frac{1}{y}\right) + \frac{z-1}{2} -\frac{iz^2}{2y} +\frac{zi}{2y} +\frac{1}{4}\right)  \frac{1}{\zeta} + O(1) \] 
\[~~~~~~~~~~~~~~~~~~~~~~~~~~~~~~~~~~~~~~~~~~~~~~~~~~~~~~~~~~~~~~~~~~~~~~~~~~~~~~~~~~~~~~~~~~~~~~~~~~~~~~~~~~~~~~~~~~~~~~~~~~~~~~~~~~~~~~~~~~~~~~~~~~~~~~~~~~~~~~~~~~~~~~~~~~~~~~~~~~~~~~~~~~~~~~~~~~~~~~~~~~~~~~~~~~~~~~~~~~~~~~~~~(R_2)\]
$\bullet$ \boxed{\textbf{\( \displaystyle \zeta = 0 \) :}} \\[1\baselineskip]
Adding the expressions $(R_1)$ and $(R_2)$, we get the residue at $\zeta=0$: 
\[\text{Res}[ \displaystyle F_n, 0\text{]}= \frac{i}{8}\left(y -\frac{1}{y}\right) + \frac{z}{2} -\frac{iz^2}{2y} +\frac{zi}{2y} -\frac{1}{4} .\]
$\bullet$ $ \boxed{  \displaystyle \zeta = \frac{ik}{N}  :}$ \\ \\
Using the \text{cot}angent formula, \[\displaystyle \text{cot}(i \zeta) = \frac{1}{i} \left(1- \frac{1}{1-e^{2 \zeta}}\right) = i\left(1- \frac{1}{1-e^{-2 \zeta}}\right),\] and by a simple calculation of the limit we obtain
\begin{align*} 
\text{Res}[ F_n, \frac{ik}{N}] &= \lim_{\zeta \to \frac{ik}{N}} \left( \zeta - \frac{ik}{N} \right) F_n(\zeta) \\
&= \frac{1}{8\pi k}\cot\left(\pi i \frac{k}{y}\right) - \frac{1}{2\pi i k} \frac{e^{-2\pi k \frac{z}{y}} e^{2\pi \frac{k}{y}}}{1 - e^{2\pi \frac{k}{y}}} \\
&= \frac{1}{8\pi i k}\left(1 - \frac{2}{1 - e^{2\pi i k \frac{z}{y}}} \right) - \frac{1}{2\pi i k} \frac{e^{-2\pi k \frac{z}{y}} e^{2\pi \frac{k}{y}}}{1 - e^{2\pi \frac{k}{y}}} .
\end{align*} 

Since \( \displaystyle \frac{1}{ k}\text{cot}(\pi i  \frac{k}{y}) \text{ is an even expression, we obtain } \) \
\begin{align*}
\displaystyle \sum _{ \substack{k=-n \\ k \ne 0} }^{n} \text{Res[} F_n, \displaystyle \frac{ik}{N}\text{]}&=  \frac{2}{8\pi i } \sum _{ k=1 }^{n} \left(\frac{1}{k}- \frac{2}{k(1-e^{2\pi   \frac{k}{y}})}\right) - \frac{1}{2\pi i}\sum _{ k=1 }^{n} 
\frac{e^{-2\pi k \frac{z}{y}} e^{2\pi \frac{k}{y}} }{k(1- e^{2\pi \frac{k}{y}})} -\frac{1}{2\pi i }\sum _{ k=1 }^{n} 
\frac{e^{2\pi k \frac{z}{y}} e^{-2\pi \frac{k}{y}} }{-k(1- e^{-2\pi \frac{k}{y}})}  \\
&= \frac{1}{4\pi i } \sum _{ k=1 }^{n} \frac{1}{k}- \frac{1}{2\pi i } \sum _{ k=1 }^{n} \frac{1}{k(1-e^{2\pi   \frac{k}{y}})} - \frac{1}{2\pi i}\sum _{ k=1 }^{n} 
\frac{e^{-2\pi k \frac{z}{y}} e^{2\pi \frac{k}{y}} }{k(1- e^{2\pi \frac{k}{y}})} -\frac{1}{2\pi i }\sum _{ k=1 }^{n} 
\frac{e^{2\pi k \frac{z}{y}}}{k(1- e^{2\pi \frac{k}{y}})} . 
\end{align*} 
$\bullet$ $\boxed{ \displaystyle  \zeta = \frac{ky}{N}  :}$ \\ \\
Similarly, we utilize the form of the \text{cot}angent again, and compute the limit to get
\begin{align*}
\text{Res}[F_n, \frac{ky}{N}] & = \lim_{\zeta \to \frac{ky}{N}} \left(\zeta - \frac{ky}{N}\right) F_n(\zeta) \\
& = -\frac{1}{8\pi k}\cot(\pi i k y) + \frac{1}{2\pi i k}\frac{e^{2\pi i k z}}{1- e^{2\pi k y}} \\
& = -\frac{1}{8\pi i k}\left(1 - \frac{2}{1 - e^{2\pi k y}}\right) + \frac{1}{2\pi i k}\frac{e^{2\pi i k z}}{1- e^{2\pi k y}}.
\end{align*}
Similarly,
\begin{align*}
\displaystyle \sum _{ \substack{k=-n \\ k \ne 0} }^{n} \text{Res[} F_n, \displaystyle \frac{ky}{N}\text{]}&=  -\frac{2}{8\pi i } \sum _{ k=1 }^{n} \left(\frac{1}{k}- \frac{2}{k(1-e^{2\pi ky}  )}\right) + \frac{1}{2\pi i}\sum _{ k=1 }^{n} 
\frac{e^{2\pi i k z}}{k(1- e^{2\pi ky})} +\frac{1}{2\pi i }\sum _{ k=1 }^{n} 
\frac{e^{-2\pi i k z}}{-k(1- e^{-2\pi ky})}  \\
&= -\frac{1}{4\pi i } \sum _{ k=1 }^{n} \frac{1}{k}+ \frac{1}{2\pi i } \sum _{ k=1 }^{n} \frac{1}{k(1-e^{2\pi ky})} + \frac{1}{2\pi i}\sum _{ k=1 }^{n} 
\frac{e^{2\pi i k z}}{k(1- e^{2\pi ky})} +\frac{1}{2\pi i }\sum _{ k=1 }^{n} 
\frac{e^{-2\pi i k z} e^{2\pi ky} }{k(1- e^{2\pi ky})} .
\end{align*} \\
Therefore, 
\begin{align*}
 2\pi i \cdot \text{Res}_{C}[F_n] & = 2\pi i  \text{ Res}[ F_n , 0 ] + 2\pi i  \sum_{\substack{k=-n \\ k \ne 0}}^{n} \text{Res}[ F_n, \frac{ik}{N} ] + \text{Res}[ F_n, \frac{ky}{N} ] \\ &= \sum_{k=1}^{n} \frac{1}{k} \frac{1}{1-e^{2k\pi y}} + \sum_{k=1}^{n} \frac{e^{2k\pi i z}}{k} \frac{1}{1-e^{2k\pi y}}  \quad + \sum_{k=1}^{n} \frac{e^{-2k\pi i z} }{k} \frac{e^{2k\pi y}}{1-e^{2k\pi y}} \\
& \quad - \sum_{k=1}^{n} \frac{1}{k} \frac{1}{1-e^{2k\pi \frac{1}{y} }} - \sum_{k=1}^{n} \frac{e^{2k\pi \frac{z}{y}}}{k} \frac{1}{1-e^{2k\pi \frac{1}{y}}}  \quad - \sum_{k=1}^{n} \frac{e^{-2k\pi \frac{z}{y}} }{k} \frac{e^{2k\pi \frac{1}{y}}}{1-e^{2k\pi \frac{1}{y}}} \\
& \quad - \frac{\pi}{4} \left(y -\frac{1}{y}\right) + \pi i z + \frac{\pi z^2}{y} - \frac{\pi z}{y} - \frac{\pi i}{2}.
\end{align*}
This completes the proof of \textit{Lemma 2}.  $ ~~~~\blacksquare$ \\ \\
\textbf{Lemma 3.} Let C be the contour in Figure \ref{fig 1}. For the sequence of functions $F_n$, we have \(
\displaystyle \lim_{n \to \infty} \zeta F_n(\zeta)= \frac{1}{8}\)  on the edges $(y,i)$, $(-y,-i),$ and  $\displaystyle \lim_{n\to \infty}\zeta F_n(\zeta)=-\frac{1}{8}$ on the other edges.\\ \\
\textit{ ~Proof.} \\ \\
We prove the \( \lim_{n \to \infty} \zeta F_n(\zeta)= \frac{1}{8}\)  on the edge $(y,i)$. \\ \\
First note that 
\begin{align} \zeta F_n(\zeta)&= -\frac{1}{8} \text{cot}(\pi i N \zeta) \text{cot}(\pi N \frac{1}{y}\zeta) +\left(\frac{1}{1-e^{2\pi N \zeta}}\right)\left(\frac{e^{-2\pi i (-\frac{zN}{y}+\frac{N}{y}i)\zeta}}{1-e^{-2\pi i \frac{N}{y}\zeta}}\right) \nonumber \\ 
&= -\left(1- \frac{1}{1-e^{-2 \pi  N \zeta}}\right) \left(1- \frac{1}{1-e^{2i \pi N \frac{1}{y}\zeta}}\right) +\left(\frac{1}{1-e^{2\pi N \zeta}}\right)\left(\frac{e^{-2\pi i (-\frac{zN}{y}+\frac{N}{y})\zeta}}{1-e^{-2\pi i \frac{N}{y}\zeta}}\right).
\end{align} 
Now the points on the edge $(y,i)$ have the form $y(1-t) +ti$ for some $0<t<1$.
Let $\zeta_0$ be a point on the edge $(y,i)$ and write $\zeta_0 = y(1-t_0) +t_0i.$ \\ \\ 
We  hence get
\begin{align*} \zeta_0 F_n(\zeta_0)&= -\frac{1}{8} \text{cot}(\pi i N \zeta_0) \text{cot}(\pi N \frac{1}{y}\zeta_0) +\left(\frac{1}{1-e^{2\pi N \zeta_0}}\right)\left(\frac{e^{-2\pi i (-\frac{zN}{y}+\frac{N}{y})\zeta_0}}{1-e^{-2\pi i \frac{N}{y}\zeta_0}}\right) \\ 
&=  \frac{1}{8} \left(1- \frac{1}{1-e^{-2 \pi  N y(1-t_0)}e^{-2 i \pi  N t_0}}\right) \left(1- \frac{1}{1-e^{2i \pi N (1-t_0)}e^{- 2 \pi N  \frac{1}{y} t_0)}}\right) \\ &~~~+\left(\frac{1}{1-e^{2\pi N y(1-t_0)}e^{2\pi i N t_0}}\right)\left(\frac{e^{-2\pi i N (-a +1 -ib) (1-t_0)}e^{2\pi \frac{N}{y} (-a+1 +ib) t_0}}{1-e^{-2\pi i N (1-t_0)}e^{2\pi  \frac{N}{y} t_0} }\right)
. 
 \end{align*} 
 However, given that $0< t_0, 1-t_0< 1$, $y>|b|\ge 0$, $N= n+\frac{1}{2}$, and $0<a, 1-a<1$, we get $\frac{1}{8}$ upon sending $n$ to $\infty.$\\
 The other limits are calculated similarly with respect to the other edges, and this completes the proof of \textit{Lemma 3}. $ ~~~~~~~~~~~~~~~~~~~~~~~~~~~~~~~~~~~~~~~~~~~~~~~~~~~~~~~~~~~~~~~~~~~~~~~~~~~~~~~~~~~~~~~~~~~~~~~~~~~~~~~~~~~~~~~~~~~~~~~~~~~~~~~~~~~~~\blacksquare$  \\ \\
We now return to the statement of the main theorem and employ \textit{Lemma 1}, \textit{lemma 2}, and \textit{Lemma 3} to complete the proof.\\ \\
\textit{Proof of theorem 1.}\\ \\
We prove  (\ref{eq:2}) for $\tau=iy$ and a fixed $z=a+ib$ where $y> |b|$ with $b<0<a<1~$ (i.e $z$ not a zero of $\vartheta_1$), then we extend the result by analytic continuation. \\
For a fixed $n$, we use the residue theorem for $F_n$ over the contour $C$ in Figure \ref{fig 1} which yields
\[
\int_{C} F_n \, d\zeta = 2\pi i \, \text{Res}_{C}[F_n].
\]
Therefore, by \textit{Lemma 2.}, we have \\
\begin{align*}
\lim_{n\to\infty}\int_{C}F_n \, d\zeta & = \lim_{n\to\infty} 2\pi i \, \text{Res}_{C}[F_n] \\
& = \sum_{m=1}^{\infty} \frac{1}{m} \frac{1}{1-e^{2m\pi y}} + \sum_{m=1}^{\infty} \frac{e^{2m\pi i z}}{m} \frac{1}{1-e^{2m\pi y}} + \sum_{m=1}^{\infty} \frac{e^{-2m\pi i z} }{m} \frac{e^{2m\pi y}}{1-e^{2m\pi y}} \\
& \quad - \sum_{m=1}^{\infty} \frac{1}{m} \frac{1}{1-e^{2m\pi \frac{1}{y} }} - \sum_{m=1}^{\infty} \frac{e^{2m\pi \frac{z}{y}}}{m} \frac{1}{1-e^{2m\pi \frac{1}{y}}}  - \sum_{m=1}^{\infty} \frac{e^{-2m\pi \frac{z}{y}} }{m} \frac{e^{2m\pi \frac{1}{y}}}{1-e^{2m\pi \frac{1}{y}}} \\
& \quad -\frac{\pi}{4}\left(y -\frac{1}{y}\right) + \pi i z +\frac{\pi z^2}{y} -\frac{\pi z}{y} -\frac{\pi i}{2} \\
& = \phi(z, iy) +\frac{\pi z^2}{y}  -\frac{\pi i}{2},
\end{align*}
where the later series converge in view of \textit{Remark 2.} \\ \\
On the other hand, computing the limit of the integral over the contour $C$ in Figure \ref{fig 1}
\[
\begin{aligned}
\lim_{n \to \infty} \int_{C} F_n(\zeta)~ d \zeta 
&= \lim_{n \to \infty} \int_{C} \zeta F_n \frac{d \zeta}{\zeta} \\
&= \int_{C}\lim_{n\to\infty}\zeta F_n \frac{d\zeta}{\zeta} \\
&= \frac{1}{8} \left( - \int_{-i}^{y} + \int_{y}^{i} -\int_{i}^{-y} +\int_{-y}^{-i}\right)\frac{d \zeta}{\zeta} \\
&= \frac{1}{4} \left( - \int_{-i}^{y} + \int_{y}^{i}\right)\frac{d \zeta}{\zeta} \\
&= \frac{1}{4} \left(-\log(y) -\frac{\pi i}{2} + \frac{\pi i}{2} -\log(y)\right) \\
&= -\log(y^{1/2}),
\end{aligned}
\]
where in the second equality we make use of the fact that $F_n$ is uniformly bounded on $C$ hence by Arzela's bounded convergence theorem we can interchange the limit with the integral, and in the third equality we invoke \textit{Lemma 3} on each edge.\\
This proves that 
\[\phi(z,iy) +\frac{\pi z^2}{y} -\frac{\pi i}{2}
= -\text{log}( y^{1/2}).\]
Now using analytic continuation with a fixed $z$ where $\vartheta_1$ does not vanish, we extend $\tau$ to the whole upper half plane, once done we extend $z$ via analytic continuation again to the whole complex plane $\mathbb{C}$. Finally, using the definition of $\phi$, we have
\[\displaystyle
\text{log}(\vartheta_1(z,\tau)) -\text{log}\left(\vartheta_1\left(\frac{z}{\tau},\frac{-1}{\tau}\right)\right) = - \text{log}(-i (-i\tau)^{1/2} e^{\pi i \frac{z^2}{\tau}}),\]
thus obtaining
\[\displaystyle
\vartheta_1\left(\frac{z}{\tau},\frac{-1}{\tau}\right)= -i (-i\tau)^{1/2} e^{\pi i \frac{z^2}{\tau}} \vartheta_1(z,\tau),\]
which completes the proof of the theorem. $ ~~~~~~~~~~~~~~~~~~~~~~~~~~~~~~~~~~~~~~~~~~~~~~~~~~~~~~~~~~~~~~~~~~~~~~~~~~~~~~~~~~~~~~~~~~~~~~~~~~~~~~~~~~~~~~~\blacksquare$  
\\
\begin{center}
Acknowledgment
\end{center} 
We would like to express our deepest gratitude to our advisor Professor Wissam Raji. We also like to thank the Department of Mathematics and the Center of Advanced Mathematical Sciences
(CAMS) at the American university of Beirut (AUB) for the guidance and support we recieved from the summer research camp (SRC). Additionally, we would like to thank the referee whose valuable comments helped us enhance the content and the structure of the paper.

Department of Mathematics, American University of Beirut, Beirut, Lebanon \\
\textit{E-mail address: mmm133@mail.aub.edu} \\
\textit{E-mail address: ays11@mail.aub.edu} \\

\end{document}